\documentclass[english]{svjour3}
\usepackage[T1]{fontenc}
\usepackage[latin9]{inputenc}
\usepackage{babel}
\usepackage{textcomp}
\usepackage{bbding}
\usepackage{url}
\usepackage{amsmath}
\usepackage{amssymb}
\usepackage{graphicx}
\usepackage[unicode=true]
 {hyperref}

\makeatletter
\RequirePackage{fix-cm}

\smartqed  % flush right qed marks, e.g. at end of proof

\usepackage{hyperref}

\usepackage{breakurl}
\usepackage{algorithm}
\usepackage{algorithmic}

\hypersetup{
  colorlinks=false,
  frenchlinks=false,
  pdfborder={0 0 0},
  naturalnames=false,
  hypertexnames=false,
  bookmarksdepth=2,
  breaklinks,
  colorlinks = true,
  allcolors = siaminlinkcolor,
  urlcolor = siamexlinkcolor,
  pdftitle={Derivative-Free Superiorization: Principle and Algorithm},
 pdfauthor={\textcopyright Yair Censor},
 pdfkeywords={Derivative-free, superiorization, constrained minimization, component-wise perturbations, proximity function, bounded perturbations, monotone proximity, proximity-target curve},
 bookmarksopen,urlcolor=blue,citecolor=blue,linkcolor=blue,pdfstartview=FitH
}

\makeatother

\begin{document}
\title{Derivative-Free Superiorization: Principle and Algorithm\thanks{Edgar Garduño would like to thank the support of DGAPA-UNAM. The work
of Yair Censor is supported by the ISF-NSFC joint research program
Grant No. 2874/19. Elias S. Helou was partially supported by CNPq
grant No. 310893/2019-4.}}
\titlerunning{\textbf{Derivative-Free Superiorization: Principle and Algorithm}}
\author{Yair Censor \and Edgar Garduño \and \linebreak{}
Elias S. Helou \and Gabor T. Herman}
\authorrunning{Yair Censor \and Edgar Garduño \and Elias S. Helou \and Gabor T.
Herman}
\institute{Yair Censor \at Department of Mathematics, University of Haifa, Mt.
Carmel, Haifa 3498838, Israel\\
\email{yair@math.haifa.ac.il}\\
Edgar Garduño \at Departamento de Ciencias de la Computación, Instituto
de Investigaciones en Matemáticas Aplicadas y en Sistemas, Universidad
Nacional Autónoma de México, Cd. Universitaria, C.P. 04510, Mexico
City, Mexico\\
\email{edgargar@ieee.org}\\
Elias S. Helou \at Instituto de Ciências Matemáticas e de Computação,
Universidade of São Paulo, São Carlos, São Paulo 13566-590, Brazil\\
\email{elias@icmc.usp.br}\\
Gabor T. Herman \at Computer Science Ph.D. Program, The Graduate
Center, City University of New York, New York, NY 10016, USA\\
\email{gabortherman@yahoo.com}}
\date{Received: date / Accepted: date / Revised: October 21, 2020}
\maketitle
\begin{abstract}
The superiorization methodology is intended to work with input data
of constrained minimization problems, that is, a target function and
a set of constraints. However, it is based on an antipodal way of
thinking to what leads to constrained minimization methods. Instead
of adapting unconstrained minimization algorithms to handling constraints,
it adapts feasibility-seeking algorithms to reduce (not necessarily
minimize) target function values. This is done by inserting target-function-reducing
perturbations into a feasibility-seeking algorithm while retaining
its feasibility-seeking ability and without paying a high computational
price. A superiorized algorithm that employs component-wise target
function reduction steps is presented. This enables derivative-free
superiorization (DFS), meaning that superiorization can be applied
to target functions that have no calculable partial derivatives or
subgradients. The numerical behavior of our derivative-free superiorization
algorithm is illustrated on a data set generated by simulating a problem
of image reconstruction from projections. We present a tool (we call
it a \textit{proximity-target curve}) for deciding which of two iterative
methods is ``better'' for solving a particular problem. The plots
of proximity-target curves of our experiments demonstrate the advantage
of the proposed derivative-free superiorization algorithm. 
\end{abstract}

\section{Introduction\label{sec:Introduction}}

\subsection{The superiorization methodology (SM)}

In many applications there exist efficient iterative algorithms for
producing \emph{constraints-compatible} solutions. Often these algorithms
are \emph{perturbation resilient} in the sense that, even if certain
kinds of changes are made at the end of each iterative step, the algorithms
still produce a constraints-compatible solution. This property is
exploited in \emph{superiorization} by using such perturbations to
steer an algorithm to an output that is as constraints-compatible
as the output of the original algorithm, but is superior (not necessarily
optimal) to it with respect to a given target function.

Superiorization has a world-view that is quite different from that
of classical constrained optimization. Both in superiorization and
in classical constrained optimization there is an assumed domain $\Omega$
and a criterion that is specified by a target function $\phi$ that
maps $\Omega$ into $\mathbb{R}$. In classical optimization it is
assumed that there is a constraints set $C$ and the task is to find
an $\boldsymbol{x}\in C$ for which $\phi(\boldsymbol{x})$ is minimal
over $C$. Two difficulties with this approach are: (1) The constraints
that arise in a practical problem may not be consistent, so $C$ could
be empty and the optimization task as stated would not have a solution.
(2) Even for nonempty $C$, iterative methods of classical constrained
optimization typically converge to a solution only in the limit and
some stopping rule is applied to terminate the process. The actual
output at that time may not be in $C$ (especially if the iterative
algorithm is initialized at a point outside $C$) and, even if it
is in $C$, it is most unlikely to be a minimizer of $\phi$ over
$C$.

Both issues are handled in the superiorization approach investigated
here by replacing the constraints set $C$ by a nonnegative real-valued
\emph{proximity} \emph{function} $\mathcal{P}r_{T}$ that indicates
how incompatible a given $\boldsymbol{x}\in\Omega$ is with specified
constraints $T$. Then the merit of an actual output $\boldsymbol{x}$
of an algorithm is represented by the smallness of the two numbers
$\mathcal{P}r_{T}(\boldsymbol{x})$ and $\phi(\boldsymbol{x})$. Roughly,
if an iterative algorithm produces an output $\boldsymbol{x}$, then
its superiorized version will produce an output $\boldsymbol{x'}$
for which $\mathcal{P}r_{T}(\boldsymbol{x'})$ is not larger than
$\mathcal{P}r_{T}(\boldsymbol{x})$, but (as in-practice demonstrated)
generally $\phi(\boldsymbol{x'})$ is smaller than $\phi(\boldsymbol{x})$.

As an example, let $\Omega=\mathbb{R}^{J}$ and consider a set $T$
of constraints of the form
\begin{equation}
\left\langle \boldsymbol{d}^{i},\boldsymbol{x}\right\rangle =h_{i},\thinspace\thinspace\thinspace i=1,2,\ldots,I,\label{eq:equalities}
\end{equation}
where $\boldsymbol{d}^{i}\in\mathbb{R}^{J}$ and $h_{i}\in\mathbb{R}$,
for all $i=1,2,\ldots,I$, and $\left\langle \cdot,\cdot\right\rangle $
is the Euclidean inner product in $\mathbb{R}^{J}$. There may or
may not be an $\boldsymbol{x}\in\mathbb{R}^{J}$ that satisfies this
set of constraints, but we can always define a proximity function
for $T$ as, for example, by
\begin{equation}
\mathcal{P}r_{T}(\boldsymbol{x}):={\displaystyle \sum_{i=1}^{I}\left(\left\langle \boldsymbol{d}^{i},\boldsymbol{x}\right\rangle -h_{i}\right)^{2}}.\label{eq:residual}
\end{equation}

There are several approaches in the literature that attempt to minimize
both competing objectives $\mathcal{P}r_{T}(\boldsymbol{x})$ and
$\phi(\boldsymbol{x})$ as a way to handle constrained minimization.
The oldest one is the penalty function approach, also useful in regularization
of inverse problems \cite{ehn00}. In that approach, the constrained
minimization problem is replaced by the unconstrained minimization
of the combination $\phi(\boldsymbol{x})+\pi\mathcal{P}r_{T}(\boldsymbol{x})$,
in which $\pi\geq0$ is a penalty parameter that governs the relative
importance of minimizing the two summands. An inherent difficulty
with this is that the penalty parameter needs to be chosen by the
user. The filter method approach \cite{filter02}, among others, was
developed to avoid this difficulty. Of course, people have also applied
multiobjective minimization with two objectives (bi-objective minimization)
to the competing objectives $\mathcal{P}r_{T}(\boldsymbol{x})$ and
$\phi(\boldsymbol{x})$. None of these approaches are close in their
underlying principles to the superiorization methodology employed
in this paper.

\subsection{Derivative-free superiorization: Expanding the boundaries of superiorization
and competing with derivative-free optimization}

Our motivating purpose in this paper is to investigate the general
applicability of derivative-free superiorization (DFS) as an alternative
to previously proposed superiorization approaches. These earlier approaches
were based on generation of nonascending vectors, for target function
reduction steps, that mostly required the ability to calculate gradients
or subgradients of the target function. Observing the body of knowledge
of derivative-free optimization (DFO), see, e.g., \cite{Conn-book-2009},
we explore a DFS algorithm and demonstrate its action numerically.

In the perturbation phase of the superiorized version of a basic algorithm
we replace the target function reduction steps that depend on gradient
or subgradient calculations by steps that use a direction search technique
which does not require any form of differentiability. Continuing the
work of \cite{CHS18}, we search the neighborhood of a current point
$\boldsymbol{x}$ for a point at which the target function exhibits
nonascent. A specific scheme for doing this is described in detail
below; particularly, in Section \ref{sec:Specific-superiorization-approac}.

While this might seem a simple technical matter, the ramifications
for practical applications of the SM are important. For example, in
intensity-modulated radiation therapy treatment planning, with photons,
protons or other particles, the normal tissue complication probability
(NTCP) is a predictor of radiobiological effects for organs at risk.
The inclusion of it, or of other biological functions, as an objective
function in the mathematical problem modeling and the planning algorithm,
is hampered because they are, in general, empirical functions whose
derivatives cannot be calculated, see, e.g., \cite{gay-niemierko-ntcp-2007}.
In the recent paper \cite{nystorm-2020} the authors list issues of
immediate clinical and practical relevance to the Proton Therapy community,
highlighting the needs for the near future but also in a longer perspective.
They say that ``...practical tools to handle the variable biological
efficiency in Proton Therapy are urgently demanded...''.

The output of a superiorized version of a constraints-compatibility-seeking
algorithm will have smaller (but not minimal) target function $\phi$
value than the output by the same constraints-compatibility-seeking
algorithm without perturbations, everything else being equal. Even
though superiorization is not an exact minimization method, we think
of it as an applicable (and possibly, more efficacious) alternative
to derivative-free constrained minimization methods applied to the
same data for two main reasons: its ability to handle constraints
and its ability to cope with very large-size problems. This is in
contrast with the current state of the art, which is as follows.

The review paper of Rios and Sahinidis \cite{Rios} ``... addresses
the solution of \emph{bound-constrained} optimization problems using
algorithms that require only the availability of objective function
values but no derivative information,'' with bound constraints imposed
on the vector $\boldsymbol{x}$. The book by Conn, Scheinberg and
Vicente \cite{Conn-book-2009} deals only with derivative-free unconstrained
minimization, except for its last chapter (of 10 pages out of the
275) entitled ``Review of constrained and other extensions to derivative-free
optimization.'' Li \emph{et al}. \cite{LCLLLL} do not even mention
constraints. In \cite{diniz2011} the numerical work deals with: ``The
dimension of the problems {[}i.e., the size of the vector $\boldsymbol{x}${]}
varies between 2 and 16, while the number of constraints are between
1 and 38, exceeding 10 in only 5 cases.'' In \cite{dfo-4-oil} the
numerical tests are limited to: ``The first case has 80 optimization
variables {[}i.e., the size of the vector $\boldsymbol{x}${]} and
only bound constraints, while the second example is a generally constrained
production optimization involving 20 optimization variables and 5
general constraints.'' Similar orders of magnitude for problem sizes
appear in the numerical results presented in \cite{Audet-Dennis-2009}
and also in the book of Audet and Hare \cite{Audet-book-2017}.

This indicates that (i) much of the literature on derivative-free
minimization is concerned with unconstrained minimization or with
bound-constraints on the variables, and (ii) many, if not all, proposed
methods were designed (or, at least, demonstrated) only for small-scale
problems. In contrast, the DFS method proposed here can handle any
type of constraints for which a separate efficient constraints-compatibility-seeking
algorithm is available and is capable of solving very large problems.
In the matter of problem sizes, we discover here, admittedly with
a very preliminary demonstration, that DFS can compete well with DFO
on large problems. Since the constraints-compatibility-seeking algorithm
forms part of the proposed DFS method, the method can use exterior
initialization (that is initializing the iterations at any point in
space). Furthermore, very large-scale problems can be accommodated.

The progressive barrier (PB) approach, described in Chapter 12 of
the book \cite{Audet-book-2017}, originally published in \cite{Audet-Dennis-2009},
is an alternative to the exterior penalty (EP) approach that we mention
in Subsection \ref{subsec:EP-approach} below. However, the PB differs
from our DFS method, in spite of some similarities with it, as we
explain in Subsection \ref{subsec:The-progressive-barrier} below.

\subsection{Earlier work on superiorization and the ``guarantee problem''}

A comprehensive overview of the state of the art and current research
on superiorization appears in our continuously updated bibliography
Internet page that currently contains 109 items \cite{sup-bib}. Research
works in this bibliography include a variety of reports ranging from
new applications to new mathematical results on the foundations of
superiorization. A special issue entitled: ``Superiorization: Theory
and Applications'' of the journal Inverse Problems \cite{Sup-Special-Issue-2017}
contains several interesting papers on the theory and practice of
SM, such as \cite{Cegielski-2017}, \cite{He2017}, \cite{hoseini-2019},
\cite{reich-zalas-2016} and \cite{Reich2017}, to name but a few.
Later papers continue research on perturbation resilience, which lies
at the heart of the SM, see, e.g., \cite{Bargetz2018}. An early paper
on bounded perturbation resilience is \cite{But06}, a recent book
containing results on the behavior of algorithms in the presence of
summable perturbations is \cite{ZAS18}.

In \cite[Section 3]{censor-levy-2019} we gave a precise definition
of the ``guarantee problem'' of the SM. We wrote there:
\begin{quotation}
The SM interlaces into a feasibility-seeking basic algorithm target
function reduction steps. These steps cause the target function to
reach lower values locally, prior to performing the next feasibility-seeking
iterations. A mathematical guarantee has not been found to date that
the overall process of the superiorized version of the basic algorithm
will not only retain its feasibility-seeking nature but also preserve
globally the target function reductions. We call this fundamental
question of the SM ``the guarantee problem of the SM''\ which is:
``under which conditions one can guarantee that a superiorized version
of a bounded perturbation resilient feasibility-seeking algorithm
converges to a feasible point that has target function value smaller
or equal to that of a point to which this algorithm would have converged
if no perturbations were applied \textendash{} everything else being
equal.''\ 
\end{quotation}
Numerous works that are cited in \cite{sup-bib} show that this global
function reduction of the SM occurs in practice in many real-world
applications. But until the guarantee problem of the SM is answered
one wonders if the SM is just a successful heuristic or if there is
a mathematical foundation for the accumulating reports on its performance
success?'' Except for a partial answer in \cite{censor-levy-2019}
with the aid of the ``concentration of measure'' principle there
are also the partial result of \cite[Theorem 4.1]{cz3-2015} about
strict Fejér monotonicity of sequences generated by an SM algorithm.

\subsection{Structure of the paper}

In Section \ref{sect:SM} we present the basics of the superiorization
methodology. We present our DFS algorithm in Section \ref{sec:Specific-superiorization-approac}
and juxtapose it with an existing superiorization algorithm that uses
derivative information. In Section \ref{sec:proximity-target-curve}
we present a tool (we call it a \textit{proximity-target curve}) for
deciding which of two iterative methods is ``better'' for solving
a particular problem. The experimental demonstration of our DFS algorithm
appears in Section \ref{sec:Experimental}. In Section \ref{sec:Discussion-and-conclusions}
we offer a brief discussion and some conclusions.

\section{\label{sect:SM}The basics of the superiorization methodology}

We follow the approach of \cite{HGDC12}. A word about terminology
before we begin: The SM has been developed with the terminology presented
here, see, e.g., \cite{Herman-JANO}. Adhering to it will assist readers
when referring to other publications and will contribute to separate
it from similar notions that are used in optimization theory.

$\Omega$ denotes a nonempty set in the Euclidean space $\mathbb{R}^{J}$.
$\mathbb{T}$ is a \textit{problem set}; each problem $T\in\mathbb{T}$
is described by a particular set of constraints such as provided,
for example, in \eqref{eq:equalities}. $\mathcal{P}r$ is a \textit{proximity
function} on $\mathbb{T}$ such that, for every $T\in\mathbb{T}$,
$\mathcal{P}r_{T}:\Omega\rightarrow\mathbb{R}_{+}$ (nonnegative real
numbers). $\mathcal{P}r_{T}\left(\boldsymbol{x}\right)$ measures
how incompatible $\boldsymbol{x}$ is with the constraints of $T$.
A \textsl{problem structure} is a pair $(\mathbb{T},\mathcal{P}r)$,
where $\mathbb{T}$ is a problem set and $\mathcal{P}r$ is a proximity
function on $\mathbb{T}$. For an $\boldsymbol{x}\in\Omega$, we say
that $\boldsymbol{x}$ is $\varepsilon$\emph{-compatible} with $T$
if $\mathcal{P}r_{T}\left(\boldsymbol{x}\right)\leq\varepsilon$.
We assume that we have computer procedures that, for any $\boldsymbol{x}\in\mathbb{R}^{J}$,
determine whether $\boldsymbol{x}\in\Omega$ and, for any $\boldsymbol{x}\in\Omega$
and $T\in\mathbb{T}$, calculate $\mathcal{P}r_{T}\left(\boldsymbol{x}\right)$.
In many applications, each problem $T\in\mathbb{T}$ is determined
by a family of sets $\left\{ C_{i}\right\} _{i=1}^{I}$, where each
$C_{i}$ is a nonempty, often closed and convex, subset of $\Omega$
and the problem $T$ is to find a point that is in the intersection
of the $C_{i}$.

We introduce $\Delta$, such that $\Omega\subseteq\Delta\subseteq\mathbb{R}^{J}$
and a \emph{target function} $\phi:\Delta\rightarrow\mathbb{R}$,
which is referred to as an optimization criterion in \cite{HGDC12}.
We assume that we have a computer procedure that, for any $\boldsymbol{x}\in\mathbb{R}^{J}$,
determines whether $\boldsymbol{x}\in\Delta$ and, if so, calculates
$\phi\left(\boldsymbol{x}\right)$.

An \textsl{algorithm} $\mathbf{P}$ \textsl{for a problem structure}
$(\mathbb{T},\mathcal{P}r)$ assigns to each problem $T\in\mathbb{T}$
a computable \emph{algorithmic operator} $\mathbf{P}_{T}:\Delta\rightarrow\Omega$.
For any \emph{initial point} $\boldsymbol{x}\in\Omega$, $\mathbf{P}_{T}$
produces the infinite sequence $\left(\left(\mathbf{P}_{T}\right)^{k}\boldsymbol{x}\right)_{k=0}^{\infty}$
of points in $\Omega$. The next definition gives a name to the first
element in a sequence $\left(\boldsymbol{x}^{k}\right)_{k=0}^{\infty}$
with $\mathcal{P}r_{T}\left(\boldsymbol{x}^{k}\right)\leq\varepsilon.$ 
\begin{definition}
\textbf{The} $\varepsilon$\textbf{-output of a sequence}

For a problem structure $(\mathbb{T},\mathcal{P}r)$, a $T\in\mathbb{T}$,
an $\varepsilon\in\mathbb{R}_{+}$ and a sequence $R:=\left(\boldsymbol{x}^{k}\right)_{k=0}^{\infty}$
of points in $\Omega$, we use $O\left(T,\varepsilon,R\right)$ to
denote the $\boldsymbol{x}\in\Omega$ that has the following properties:
$\mathcal{P}r_{T}(\boldsymbol{x})\leq\varepsilon,$ and there is a
nonnegative integer $K$ such that $\boldsymbol{x}^{K}=\boldsymbol{x}$
and, for all nonnegative integers $k<K$, $\mathcal{P}r_{T}\left(\boldsymbol{x}^{k}\right)>\varepsilon$.
If there is such an $\boldsymbol{x}$, then it is unique. If there
is no such $\boldsymbol{x}$, then we say that $O\left(T,\varepsilon,R\right)$
is \textit{undefined}, otherwise it is \textit{defined}.
\end{definition}
If $R$ is an infinite sequence generated by a process that repeatedly
applies $\mathbf{P}_{T}$, then $O\left(T,\varepsilon,R\right)$ is
the \textit{output} produced by that process when we add to it instructions
that make it terminate as soon as it reaches a point that is $\varepsilon$-compatible
with $T$. Roughly, we refer to $\mathbf{P}$ as a \textit{feasibility-seeking
algorithm} for a problem structure $(\mathbb{T},\mathcal{P}r)$ that
arose from a particular application if, for all $T\in\mathbb{T}$
and $\varepsilon\in\mathbb{R}_{+}$ of interest for the application,
$O\left(T,\varepsilon,R\right)$ is defined for all infinite sequences
$R$ generated by repeated applications $\mathbf{P}_{T}$. Each application
of $\mathbf{P}_{T}$ is referred to as a \textit{feasibility-seeking
step}.
\begin{definition}
\textbf{Strong perturbation resilience}
\end{definition}
An algorithm $\mathbf{P}$ for a problem structure $(\mathbb{T},\mathcal{P}r)$
is said to be \textit{strongly perturbation resilient} if, for all
$T\in\mathbb{T}$,
\begin{enumerate}
\item there is an $\varepsilon\in\mathbb{R}_{+}$ such that $O\left(T,\varepsilon,\left(\left(\mathbf{P}_{T}\right)^{k}\boldsymbol{x}\right)_{k=0}^{\infty}\right)$
is defined for every $\boldsymbol{x}\in\Omega$;
\item for all $\varepsilon\in\mathbb{R}_{+}$ such that $O\left(T,\varepsilon,\left(\left(\mathbf{P}_{T}\right)^{k}\boldsymbol{x}\right)_{k=0}^{\infty}\right)$
is defined for every $\boldsymbol{x}\in\Omega$, we also have that
$O\left(T,\varepsilon',R\right)$ is defined for every $\varepsilon'>\varepsilon$
and for every sequence $R=\left(\boldsymbol{x}^{k}\right)_{k=0}^{\infty}$
of points in $\Omega$ generated by
\begin{equation}
\boldsymbol{x}^{k+1}=\mathbf{P}_{T}\left(\boldsymbol{x}^{k}+\beta_{k}\boldsymbol{v}^{k}\right),\:\mathrm{for\:all\:}k\geq0,\label{eq:perturbations}
\end{equation}
where $\beta_{k}\boldsymbol{v}^{k}$ are \textit{bounded perturbations},
meaning that the sequence $\left(\beta_{k}\right)_{k=0}^{\infty}$
of nonnegative real numbers is \textit{summable} (that is, ${\displaystyle \sum\limits _{k=0}^{\infty}}\beta_{k}\,<\infty$),
the sequence $\left(\boldsymbol{v}^{k}\right)_{k=0}^{\infty}$ of
vectors in $\mathbb{R}^{J}$ is bounded and, for all $k\geq0$, $\boldsymbol{x}^{k}+\beta_{k}\boldsymbol{v}^{k}\in\Delta$.
\end{enumerate}
Sufficient conditions for strong perturbation resilience appeared
in \cite[Theorem 1]{HGDC12}. Here and elsewhere, while the $\beta_{k}$s
are nonnegative real numbers, one should disallow the trivial case
in which all $\beta_{k}$s are zero because obviously that would completely
nullify all perturbations. Notice that for most feasibility-seeking
algorithms we have $\lim_{k\to\infty}\mathcal{P}r\left(\left(\mathbf{P}_{T}\right)^{k}\boldsymbol{x}\right)=\inf_{\boldsymbol{y}}\mathcal{P}r(\boldsymbol{y})$,
which means that superiorization of a strongly perturbation resilient
feasibility-seeking algorithm will generate an infinite perturbed
sequence such that $\lim_{k\to\infty}\mathcal{P}r\left(\boldsymbol{x}^{k}\right)=\inf_{\boldsymbol{x}}\mathcal{P}r(\boldsymbol{x})$.

With respect to the target function $\phi:\Delta\rightarrow\mathbb{R}$,
we adopt the convention that a point in $\Delta$ for which the value
of $\phi$ is smaller is considered \textit{superior} to a point in
$\Delta$ for which the value of $\phi$ is larger. The essential
idea of the SM is to make use of the perturbations of \eqref{eq:perturbations}
to transform a strongly perturbation resilient algorithm that seeks
a constraints-compatible solution (referred to as the \emph{Basic
Algorithm}) into a \emph{superiorized version} whose outputs are equally
good from the point of view of constraints-compatibility, but are
superior (not necessarily optimal) with respect to the target function
$\phi$. This can be done by making use of the following concept.
\begin{definition}
\label{def:old-nonascent}\cite{HGDC12} \textbf{Nonascending vector}
\end{definition}
Given a function $\phi:\Delta\rightarrow\mathbb{R}$ and a point $\boldsymbol{y}\in\mathbb{R}^{J}$,
we say that a $\boldsymbol{d}\in\mathbb{R}^{J}$ is a \emph{nonascending
vector for} $\phi$ \emph{at} $\boldsymbol{y}$ if $\left\Vert \boldsymbol{d}\right\Vert \leq1$
and there is a $\delta>0$ such that 
\begin{equation}
\text{for all }\lambda\in\left[0,\delta\right]\text{ we have }\phi\left(\boldsymbol{y}+\lambda\boldsymbol{d}\right)\leq\phi\left(\boldsymbol{y}\right).\label{eq:nonascend-1}
\end{equation}

Obviously, the zero vector $\mathbf{0}$ (all components are 0) is
always such a vector, but for the SM to work we need a strict inequality
to occur in \eqref{eq:nonascend-1} frequently enough. Generation
of nonascending vectors, used for target function reduction steps,
has been based mostly on the following theorem or its variants such
as \cite[Theorem 1]{GH14} and \cite[unnumbered Theorem on page 7]{GTH},
which provide sufficient conditions for a nonascending vector.
\begin{theorem}
\label{thm:old-thm}\cite[Theorem 2]{HGDC12}. Let $\phi:\mathbb{R}^{J}\rightarrow\mathbb{R}$
be a convex function and let $\boldsymbol{x}\in\mathbb{R}^{J}$. Let
$\boldsymbol{g}\in\mathbb{R}^{J}$ satisfy the property: For $1\leq j\leq J$,
if the $j$th component $g_{j}$ of $\boldsymbol{g}$ is not zero,
then the partial derivative $\frac{\partial\phi}{\partial x_{j}}(\boldsymbol{x})$
of $\phi$ at $\boldsymbol{x}$ exists and its value is $g_{j}$.
Define $\boldsymbol{d}$ to be the zero vector if $\left\Vert \boldsymbol{g}\right\Vert =0$
and to be $-\boldsymbol{g}/\left\Vert \boldsymbol{g}\right\Vert $
otherwise. Then $\boldsymbol{d}$ is a nonascending vector for $\phi$
at $\boldsymbol{x}$.
\end{theorem}
In order to use this theorem, $\phi$ must have at least one calculable
partial derivative (which is nonzero) at points in the domain of $\phi$.
Otherwise, the theorem would apply only to the zero vector, which
is a useless nonascending vector because it renders the SM ineffective.
If $\phi$ is not differentiable at some points we can just take $\boldsymbol{v}^{k}=0$
at those points, but $\phi$ needs to be differentiable at ``enough''
points for the superiorization to be effective. To enable application
of the SM to target functions that have no calculable partial derivatives
or subgradients, we proposed in \cite{CHS18} to search for a point
in the neighborhood of $\boldsymbol{x}$ at which the target function
exhibits nonascent by comparing function values at points of a fixed
distance from $\boldsymbol{x}$ along the space coordinates. To obtain
a sequence of nonascending points without making use of Theorem \ref{thm:old-thm},
we replaced in \cite{CHS18} the notion of a nonascending vector by
the following alternative notion.
\begin{definition}
\label{def:nonascent}\cite{CHS18} \textbf{Nonascending $\delta$-bound
direction}
\end{definition}
Given a target function $\phi:\Delta\rightarrow\mathbb{R}$ where
$\Delta\subseteq\mathbb{R}^{J}$, a point $\boldsymbol{y}\in\Delta$,
and a positive $\delta\in\mathbb{R}$, we say that $\boldsymbol{d}\in\mathbb{R}^{J}$
is a \emph{nonascending} $\delta$-\emph{bound direction for} $\phi$
\emph{at} $\boldsymbol{y}$ if $\|\boldsymbol{d}\|\leq\delta$, $\boldsymbol{y}+\boldsymbol{d}\in\Delta$
and $\phi(\boldsymbol{y}+\boldsymbol{d})\leq\phi(\boldsymbol{y})$.
The collection of all such vectors is called a \emph{nonascending}
$\delta$-\emph{ball} and is denoted by $\mathcal{B}_{\delta,\phi}(\boldsymbol{y})$,
that is, 
\begin{equation}
\mathcal{B}_{\delta,\phi}(\boldsymbol{y}):=\left\{ \boldsymbol{d}\in\mathbb{R}^{J}\mid\|\boldsymbol{d}\|\leq\delta,\;(\boldsymbol{y}+\boldsymbol{d})\in\Delta,\ \phi(\boldsymbol{y}+\boldsymbol{d})\leq\phi(\boldsymbol{y})\right\} .\label{eq:ball}
\end{equation}
The zero vector is contained in each nonascending $\delta$-ball,
that is, $\boldsymbol{0}\in\mathcal{B}_{\delta,\phi}(\boldsymbol{y})$
for each $\delta>0$ and $\boldsymbol{y}\in\Delta$. The purpose of
this definition is to allow the use, as a direction of target function
decrease, of any vector $\boldsymbol{d}\in\mathbb{R}^{J}$ for which
$\phi(\boldsymbol{y}+\boldsymbol{d})\leq\phi(\boldsymbol{y})$ holds
locally only for $\boldsymbol{d}$, and not throughout a certain interval
as in Definition \ref{def:old-nonascent}. The vector $\boldsymbol{d}$
depends on the value of $\delta$ and they may be determined simultaneously
in the superiorization process, as seen below. This kind of nonascent
was referred to as \emph{local nonascent} in \cite[Subsection 2.3]{CHS18}.
Obviously, local nonascent is a more general notion since every nonascending
vector according to Definition \ref{def:old-nonascent} is also a
nonascending $\delta$-\textit{bound direction} according to Definition
\ref{def:nonascent} but not vice versa. The advantage of this notion
is that it is detectable by using only function value calculations. 

The following easily-proved proposition unifies these approaches in
the convex case. 
\begin{proposition}
Let $\phi:\mathbb{R}^{J}\rightarrow\mathbb{R}$ be a convex function
and let $\boldsymbol{x}\in\mathbb{R}^{J}$. If $\boldsymbol{d}\in\mathbb{R}^{J}$
is a nonascending $\delta$-bound direction for $\phi$ at $\mathbf{\boldsymbol{x}}$,
then either $\boldsymbol{d}=\boldsymbol{0}$ (and hence $\boldsymbol{d}$
is a nonascending vector for $\phi$ at $\boldsymbol{x}$) or $\boldsymbol{d}/\|\boldsymbol{d}\|$
is a nonascending vector for $\phi$ at $\boldsymbol{x}$.
\end{proposition}
The idea of calculating $\delta$ (equivalently, the step-size $\gamma_{\ell}$
in the superiorized algorithms presented in the next section) simultaneously
with a direction of nonascent appeared in a completely different way
in \cite{LZZS16}, where they use an additional internal loop of a
penalized minimization to calculate the direction of nonascent; see
also \cite{LZZSW19}.

\section{\label{sec:Specific-superiorization-approac}Specific superiorization
approaches}

This section presents two specific approaches to superiorizing a Basic
Algorithm that operates by repeated applications of an algorithmic
operator $\mathbf{P}_{T}$ starting from some initial point. The first
approach produces the superiorized version that is named \textbf{Algorithm
\ref{alg:NonAsc_Superiorization}} below, it has been published in
the literature previously \cite[page 5537]{HGDC12}. The second approach,
named \textbf{Algorithm \ref{alg:Component-wise}} below, is novel
to this paper.

The two superiorized versions have some things in common. They are
both iterative procedures in which $k$ is used as the iteration index.
The first two steps of both algorithms sets $k$ to $0$ and $\boldsymbol{x}^{0}$
to a given initial vector $\boldsymbol{\bar{x}}\in\Delta$. They both
assume that we have available a summable sequence $\left(\gamma_{\ell}\right)_{\ell=0}^{\infty}$
of nonnegative real numbers, not all of them zero, (for example, $\gamma_{\ell}=a^{\ell}$,
where $0<a<1$). In Step 3 of both algorithms, $\ell$ is initialized
to $-1$ (this is acceptable since $\ell$ is increased by 1 before
the first time $\gamma_{\ell}$ is used). In both algorithms the iterative
step that produces $\boldsymbol{x}^{k+1}$ from $\boldsymbol{x}^{k}$,
as in \eqref{eq:perturbations}, is specified within a \textbf{repeat}
loop that first performs a user-specified number, $N$, of \emph{perturbation
steps} followed by one \emph{feasibility-seeking step} that uses the
algorithmic operator $\mathbf{P}_{T}$. In more detail, the \textbf{repeat}
loop in each of the algorithms has the following form. After initializing
the loop index $n$ to $0$ and setting $\boldsymbol{x}^{k,0}$ to
$\boldsymbol{x}^{k}$, it produces one-by-one $\boldsymbol{x}^{k,1},\boldsymbol{x}^{k,2}$,
$\:\ldots,\:\boldsymbol{x}^{k,N}$ (these are the iterations of the
perturbation steps), followed by producing $\boldsymbol{x}^{k+1}=\mathbf{P}_{T}\boldsymbol{x}^{k,N}$
(the feasibility-seeking step). The difference between the two algorithms
is in how they perform the perturbations for getting from $\boldsymbol{x}^{k,n}$
to $\boldsymbol{x}^{k,n+1}$.

\begin{algorithm}

\caption{Superiorization using nonascending vectors}

\label{alg:NonAsc_Superiorization}

\begin{algorithmic}[1]

\STATE{\textbf{set} $k=0$}

\STATE{\textbf{set} $\boldsymbol{x}^{k}=\boldsymbol{\bar{x}}$}

\STATE{\textbf{set} $\ell=-1$}

\STATE{\textbf{repeat}}

\STATE{$\quad$\textbf{set} $n=0$}

\STATE{$\quad$\textbf{set} $\boldsymbol{x}^{k,n}=\boldsymbol{x}^{k}$}

\STATE{$\quad$\textbf{while} $n<N$}

\STATE{$\quad\quad$\textbf{set} $\boldsymbol{v}^{k,n}$ to be a
nonascending vector for $\phi$ at $\boldsymbol{x}^{k,n}$}

\STATE{$\quad\quad$\textbf{set} \emph{loop}=\emph{true}}

\STATE{$\quad\quad$\textbf{while} \emph{loop}}

\STATE{$\quad\quad\quad$\textbf{set} $\ell=\ell+1$}

\STATE{$\quad\quad\quad$\textbf{set} $\boldsymbol{z}=\boldsymbol{x}^{k,n}+\gamma_{\ell}\boldsymbol{v}^{k,n}$}

\STATE{$\quad\quad\quad$\textbf{if} $\boldsymbol{z}\in\Delta$ \textbf{and}
$\phi\left(\boldsymbol{z}\right)\leq\phi\left(\boldsymbol{x}^{k}\right)$
\textbf{then}}

\STATE{$\quad\quad\quad\quad$\textbf{set} $\boldsymbol{x}^{k,n+1}=\boldsymbol{z}$}

\STATE{$\quad\quad\quad\quad$\textbf{set} $n=n+1$}

\STATE{$\quad\quad\quad\quad$\textbf{set} \emph{loop} = \emph{false}}

\STATE{$\quad$\textbf{set} $\boldsymbol{x}^{k+1}=\mathbf{P}_{T}\boldsymbol{x}^{k,N}$}

\STATE{$\quad$\textbf{set} $k=k+1$}

\end{algorithmic}

\end{algorithm}

We state an important property of \textbf{Algorithm \ref{alg:NonAsc_Superiorization}};
for a proof see \cite[Section II.E]{HGDC12}.
\begin{theorem}
\label{perturbatuion_resilience}Suppose that the algorithm $\mathbf{P}$
for a problem structure $(\mathbb{T},\mathcal{P}r)$ is strongly perturbation
resilient. Suppose further that $T\in\mathbb{T}$ and $\varepsilon\in\mathbb{R}_{+}$
are such that $O\left(T,\varepsilon,\left(\left(\mathbf{P}_{T}\right)^{k}\boldsymbol{x}\right)_{k=0}^{\infty}\right)$
is defined for every $\boldsymbol{x}\in\Omega$. It is then the case
that $O\left(T,\varepsilon',R\right)$ is defined for every $\varepsilon'>\varepsilon$
and every sequence $R=\left(\boldsymbol{x}^{k}\right)_{k=0}^{\infty}$
produced by \textbf{Algorithm \ref{alg:NonAsc_Superiorization}.}
\end{theorem}
The pseudo-code of \textbf{Algorithm \ref{alg:NonAsc_Superiorization}}
does not specify how the nonascending vector in Step 8 is to be selected.
In publications using \textbf{Algorithm \ref{alg:NonAsc_Superiorization}},
such details are usually based on a variant of Theorem \eqref{thm:old-thm},
resulting in a not derivative-free algorithm.

For the specification of \textbf{Algorithm \ref{alg:Component-wise}}
we let, for $1\leq j\leq J$, $\mathbf{e}^{j}$ be the vector in $\mathbb{R}^{J}$
all of whose components are $0$, except for the $j$th component,
which is $1$. The set of \textit{coordinate directions} is defined
as $\Gamma:=\left\{ \mathbf{e}^{j}\:|\:1\leq j\leq J\right\} \,\cup\:\left\{ -\mathbf{e}^{j}\:|\:1\leq j\leq J\right\} $.
We assume that $\left(\mathbf{c}^{m}\right)_{m=0}^{\infty}$ is a
given sequence of coordinate directions such that any subsequence
of length $2J$ contains $\Gamma$.

The component-wise approach indicated in \textbf{Algorithm \ref{alg:Component-wise}}
has been presented previously in the literature as a stand-alone optimization
method, referred to as ``compass search''; see the review paper
by Kolda \textit{et al}. \cite{Kolda2003}, which is an excellent
source for direct search methods that do not explicitly use derivatives.
Here we apply it to superiorize a feasibility-seeking algorithm, rather
than as a stand-alone optimization method.

\begin{algorithm}

\caption{Component-wise derivative-free superiorization}

\label{alg:Component-wise}

\begin{algorithmic}[1]

\STATE{\textbf{set} $k=0$}

\STATE{\textbf{set} $\boldsymbol{x}^{k}=\boldsymbol{\bar{x}}$}

\STATE{\textbf{set} $\ell=-1$}

\STATE{\textbf{set} $m=-1$}

\STATE{\textbf{repeat}}

\STATE{$\quad$\textbf{set} $n=0$}

\STATE{$\quad$\textbf{set} $\boldsymbol{x}^{k,n}=\boldsymbol{x}^{k}$}

\STATE{$\quad$\textbf{while} $n<N$}

\STATE{$\quad\quad$\textbf{set} $\boldsymbol{x}^{k,n+1}=\boldsymbol{x}^{k,n}$}

\STATE{$\quad\quad$\textbf{set} $\ell=\ell+1$}

\STATE{$\quad\quad$\textbf{set} $L=-1$}

\STATE{$\quad\quad$\textbf{while} $L<2J$}

\STATE{$\quad\quad\quad$\textbf{set} $L=L+1$}

\STATE{$\quad\quad\quad$\textbf{set} $m=m+1$}

\STATE{$\quad\quad\quad$\textbf{set} $\boldsymbol{z}=\boldsymbol{x}^{k,n}+\gamma_{\ell}\mathbf{c}^{m}$}

\STATE{$\quad\quad\quad$\textbf{if} $\boldsymbol{z}\in\Delta$ \textbf{and}
$\phi\left(\boldsymbol{z}\right)<\phi\left(\boldsymbol{x}^{k,n}\right)$
\textbf{then}}

\STATE{$\quad\quad\quad\quad$\textbf{set} $\boldsymbol{x}^{k,n+1}=\boldsymbol{z}$}

\STATE{$\quad\quad\quad\quad$\textbf{set} $L=2J$}

\STATE{$\quad\quad$\textbf{set} $n=n+1$}

\STATE{$\quad$\textbf{set} $\boldsymbol{x}^{k+1}=\mathbf{P}_{T}\boldsymbol{x}^{k,N}$}

\STATE{$\quad$\textbf{set} $k=k+1$}

\end{algorithmic}

\end{algorithm}

We make the following comments:
\begin{enumerate}
\item Steps 15, 16 and 17 of \textbf{Algorithm \ref{alg:Component-wise}}
implement nonascending $\gamma_{\ell}$-bound directions, as in Definition
\ref{def:nonascent}. In doing so, \textbf{Algorithm \ref{alg:Component-wise}}
realizes in a component-wise manner the algorithmic framework of \cite{CHS18}
(specifically, as expressed in Steps 7 and 8 of Algorithm 1 in that
paper).
\item No partial derivatives are used by \textbf{Algorithm \ref{alg:Component-wise}}.
\item Step 16 of \textbf{Algorithm \ref{alg:Component-wise}} is similar
to Step 13 of \textbf{Algorithm \ref{alg:NonAsc_Superiorization}}.
One difference is the use of strict inequality in \textbf{Algorithm
\ref{alg:Component-wise}}, the reason for this is that it was found
advantageous in some applications of the algorithm. In addition, the
\textbf{while} loop due to Step 12 of \textbf{Algorithm \ref{alg:Component-wise}}
is executed at most $2J$ times, but there is no upper bound on the
(known to be finite) number of executions of the \textbf{while} loop
due to Step 10 of \textbf{Algorithm \ref{alg:NonAsc_Superiorization}}.
Also, it follows from the pseudo-code of \textbf{Algorithm \ref{alg:Component-wise}}
that, for all $k\geq0$ and $0\leq n\leq N$, $\phi\left(\boldsymbol{x}^{k,n}\right)<\phi\left(\boldsymbol{x}^{k}\right)$,
even though there is no explicit check for this as in Step 13 of \textbf{Algorithm
\ref{alg:NonAsc_Superiorization}}.
\item \textbf{Algorithm \ref{alg:Component-wise}} has the essential property
that it cannot get stuck in a particular iteration $k$ because the
value of $L$ increases in an execution of the \textbf{while} loop
of Step 13 and the value of $n$ increases in an execution of the
\textbf{while} loop of Step 8.
\item \textbf{Algorithm \ref{alg:Component-wise}} shares with \textbf{Algorithm
\ref{alg:NonAsc_Superiorization}} the important property in Theorem
\ref{perturbatuion_resilience}. Stated less formally: ``For a strongly
perturbation resilient algorithm, if for all initial points from $\Omega$
the infinite sequence produced by an algorithm contains an $\varepsilon$-compatible
point, then all perturbed sequences produced by the superiorized version
of the algorithm contain an $\varepsilon'$-compatible point, for
any $\varepsilon'>\varepsilon$.''
\item At present there is no mathematical proof to guarantee that the output
of a superiorized version of a constraints-compatibility-seeking algorithm
will have smaller target function $\phi$ value than the output by
the same constraints-compatibility-seeking algorithm without perturbations,
everything else being equal. A partial mathematical result toward
coping with this lacuna, in the framework of weak superiorization,
is provided by Theorem 4.1 in \cite{cz3-2015}.\footnote{The approach followed in the present paper was termed \textsl{strong
superiorization} in \cite[Section 6]{cz3-2015} and \cite{Cen15}
to distinguish it from \textsl{weak superiorization}, wherein asymptotic
convergence to a point in $C$ is studied instead of $\varepsilon$-compatibility.}
\end{enumerate}

\section{\label{sec:proximity-target-curve}The proximity-target curve}

We now give a tool for deciding which of two iterative methods is
``better'' for solving a particular problem. Since an iterative
method produces a sequence of points, our definition is based on such
sequences. Furthermore, since we are interested in the values of two
functions (proximity function and target function) at each of the
points, the efficacy of the behavior of the iterative method can be
represented by a curve in two-dimensional space, defined as the proximity
target curve below. It indicates the target value for any achieved
proximity value. This leads to the intuitive concept of an algorithm
being ``better'' than another one, if its proximity target curve
is below that of the other one (that is, the target value for it is
always smaller than the target value of the other one for the same
proximity value). Such may not always be the case, the two proximity
curves may cross each other, providing us with intervals of proximity
values within which one or the other method is better.

This way of thinking is commonly used in many fields of science in
situations where it is desirable to obtain an object for which the
values of two evaluating functions are simultaneously small. A prime
example is in estimation theory where we desire an estimation method
with both small bias and small variance. More specifically, and nearer
to the application areas of the authors, is the concept of a receiver
operating characteristics (ROC) curve that illustrates the diagnostic
ability of a binary classifier system as its discrimination threshold
is varied. It is created by plotting the true positive rate against
the false positive rate at various threshold settings. One classifier
system is considered ``better'' than the other one if its ROC curve
is above that of the other one; but, just as for our proximity-target
curves, the ROC curves for two classifier systems may cross each other.
There are many publications on the role of ROC curves in the evaluation
of medical imaging techniques; see, for example, \cite{METZ89a,SWET79a}.
Their use for image reconstruction algorithm evaluation is discussed,
for example, in \cite{COOL92a}.

For incarnations of the definitions given in this section, the reader
may wish to look ahead to Figure \ref{fig:Target-function-values}.
That figure illustrates the notions discussed in this section for
two particular finite sequences $R:=\left(\boldsymbol{x}^{k}\right)_{k=K_{lo}}^{K_{hi}}$
and $S:=\left(\boldsymbol{y}^{k}\right)_{k=L_{lo}}^{L_{hi}}$. The
details of how those sequences were produced are given below in Subsection
\ref{subsec:Algorithmic-details-and}.
\begin{definition}
\textbf{Monotone proximity of a finite sequence}
\end{definition}
For a problem structure $(\mathbb{T},\mathcal{P}r)$, a $T\in\mathbb{T}$,
positive integers $K_{lo}$ and $K_{hi}>K_{lo}$, the finite sequence
$R:=\left(\boldsymbol{x}^{k}\right)_{k=K_{lo}}^{K_{hi}}$ of points
in $\Omega$ is said to be of \emph{monotone proximity} if for $K_{lo}<k\leq K_{hi}$,
$\mathcal{P}r_{T}\left(\boldsymbol{x}^{k-1}\right)>\mathcal{P}r_{T}\left(\boldsymbol{x}^{k}\right)$.
\begin{definition}
\textbf{The proximity-target curve of a finite sequence}
\end{definition}
For a problem structure $(\mathbb{T},\mathcal{P}r)$, a $T\in\mathbb{T}$,
a target function $\phi:\Omega\rightarrow\mathbb{R}$, positive integers
$K_{lo}$ and $K_{hi}>K_{lo}$, let $R:=\left(\boldsymbol{x}^{k}\right)_{k=K_{lo}}^{K_{hi}}$
be a sequence of monotone proximity. Then the \emph{proximity-target
curve} $P\subseteq\mathfrak{\mathbb{R}}^{2}$ associated with $R$
is uniquely defined by:
\begin{enumerate}
\item For $K_{lo}\leq k\leq K_{hi}$, $\left(\mathcal{P}r_{T}\left(\boldsymbol{x}^{k}\right)\!,\phi\left(\boldsymbol{x}^{k}\right)\right)\in P$.
\item The intersection $\{(y,x)\in\mathbb{R}^{2}:\mathcal{P}r_{T}\left(\boldsymbol{x}^{k}\right)\leq y\leq\mathcal{P}r_{T}\left(\boldsymbol{x}^{k-1}\right)\}\cap P$
is the line segment from $\left(\mathcal{P}r_{T}\left(\boldsymbol{x}^{k-1}\right),\phi\left(\boldsymbol{x}^{k-1}\right)\right)$
to $\left(\mathcal{P}r_{T}\left(\boldsymbol{x}^{k}\right),\phi\left(\boldsymbol{x}^{k}\right)\right)$.
\end{enumerate}
\begin{definition}
\label{def:comparison}\textbf{Comparison of proximity-target curves
of sequences}

For a problem structure $(\mathbb{T},\mathcal{P}r)$, a $T\in\mathbb{T}$,
a target function $\phi:\Omega\rightarrow\mathbb{R}$, positive integers
$K_{lo}$, $K_{hi}>K_{lo}$, $L_{lo}$, $L_{hi}>L_{lo}$, let $R:=\left(\boldsymbol{x}^{k}\right)_{k=K_{lo}}^{K_{hi}}$
and $S:=\left(\boldsymbol{y}^{k}\right)_{k=L_{lo}}^{L_{hi}}$ be sequences
of points in $\Omega$ of monotone proximity for which $P$ and $Q$
are their respective associated proximity-target curves. Define
\begin{equation}
\begin{array}{c}
t:=\max\left(\mathcal{P}r_{T}\left(\boldsymbol{x}^{K_{hi}}\right),\mathcal{P}r_{T}\left(\boldsymbol{y}^{L_{hi}}\right)\right),\\
u:=\min\left(\mathcal{P}r_{T}\left(\boldsymbol{x}^{K_{lo}}\right),\mathcal{P}r_{T}\left(\boldsymbol{y}^{L_{lo}}\right)\right).
\end{array}\label{eq:limits-1}
\end{equation}
Then $R$ is \textit{better targeted} than $S$ if:
\end{definition}
\begin{enumerate}
\item $t\leq u$ and
\item for any real number $h$, if $t\leq h\leq u$, $\left(h,v\right)\in P$
and $\left(h,w\right)\in Q$, then $v\leq w$.
\end{enumerate}
\medskip{}
\quad{}Let us see how this last definition translates into something
that is intuitively desirable. Suppose that we have an iterative algorithm
that produces a sequence, $\boldsymbol{y}^{0},\boldsymbol{y}^{1},\boldsymbol{y}^{2},\cdots$,
of which $S:=\left(\boldsymbol{y}^{k}\right)_{k=L_{lo}}^{L_{hi}}$
is a subsequence. An alternative algorithm that produces a sequence
of points of which $R:=\left(\boldsymbol{x}^{k}\right)_{k=K_{lo}}^{K_{hi}}$
is a subsequence that is better targeted than $S$ has a desirable
property: Within the range $\left[t,u\right]$ of proximity values,
the point that is produced by the alternative algorithm with that
proximity value, is likely to have lower (and definitely not higher)
value of the target function as the point with that proximity value
that is produced by the original algorithm. This property is stronger
than what we stated before, namely that superiorization produces an
output that is equally good from the point of view of proximity, but
is superior with respect to the target function. Here the single output
determined by a fixed $\varepsilon$ is replaced by a set of potential
outputs for any $\varepsilon\in\left[t,u\right]$.

\section{\label{sec:Experimental}Experimental demonstration of derivative-free
component-wise superiorization}

\subsection{\label{subsec:Goal}Goal and general methodology}

Our goal is to demonstrate that component-wise superiorization (\textbf{Algorithm
\ref{alg:Component-wise}}) is a viable efficient DFS method to handle
data of constrained-minimization problems (that is, a target function
and a set of constraints), when the target function has no calculable
partial derivatives.

To ensure the meaningfulness and worthiness of our experiments, we
generate the constraints and choose a target function, that has no
calculable partial derivatives, inspired by an application area of
constrained optimization, namely image reconstruction from projections
in computerized tomography (CT). 

For the so-obtained data we consider two runs of \textbf{Algorithm
\ref{alg:Component-wise}}, one with and the other without the component-wise
perturbation steps. To be exact, by ``without perturbation'' we
mean that Steps 10\textendash 18 in \textbf{Algorithm \ref{alg:Component-wise}}
are deleted so that $\boldsymbol{x}^{k,N}=\boldsymbol{x}^{k}$, which
amounts to running the feasibility-seeking basic algorithm $\mathbf{P}_{T}$
without any perturbations. Everything else is equal in the two runs,
such as the initialization point $\boldsymbol{\bar{x}}$ and all parameters
associated with the application of the feasibility-seeking basic algorithm
in Step 20. The results are presented below by plots of proximity-target
curves that show that the target function values of \textbf{Algorithm
\ref{alg:Component-wise}} when run ``with perturbations'' are systematically
lower than those of the same algorithm without the component-wise
perturbations.

The numerical behavior of \textbf{Algorithm \ref{alg:Component-wise}},
as demonstrated by our experiment, makes it a meritorious choice for
superiorization in situations involving a derivative-free target function
and a set of constraints.

To reach the goal described above we proceed in the following stages.
\begin{enumerate}
\item Specification of a problem structure $(\mathbb{T},\mathcal{P}r)$
for the experimental demonstration, and generation of constraints,
simulated from the application of image reconstruction from projections
in computerized tomography.
\item Choice of a $\Delta$ and a derivative-free target function $\phi$
for the experiment.
\item Specification of the algorithmic operator $\mathbf{P}_{T}$ to be
used in \textbf{Algorithm \ref{alg:Component-wise}}. This is chosen
so that the Basic Algorithm that operates by repeated applications
of $\mathbf{P}_{T}$ is a standard sequential iterative projections
method for feasibility-seeking of systems of linear equations; a version
of the Algebraic Reconstruction Techniques (ART) \cite[Chapter 11]{GTH-book}
that is equivalent to Kaczmarz's projections method \cite{kacmarcz}.
\item Specification of algorithmic details and parameters, such as $N$
and $\gamma_{\ell}$ in \textbf{Algorithm \ref{alg:Component-wise}}.
\end{enumerate}

\subsection{\label{subsec:The-selected-problem}Problem selection, constraints
generation and choices of $\Delta$ and of the target function}

We generate the constraints and chose a target function from the application
area of image reconstruction from projections in computerized tomography
(CT).\footnote{The term projection has in this field a different meaning than in
convex analysis. It stands for a set of estimated line integrals through
the image that has to be reconstructed, see \cite[page 3]{GTH-book}.} The problem structure $(\mathbb{T},\mathcal{P}r)$ for our demonstration
has been used in the literature for comparative evaluations of various
algorithms for CT \cite{GTH,GTH-book,HGDC12,NDH12}. It is of the
type described in Section \ref{sec:Introduction} by \eqref{eq:equalities}
and \eqref{eq:residual}. Specifically, vectors $\boldsymbol{x}$
in $\Omega=\mathbb{R}^{J}$ represent two-dimensional (2D) images,
with each component of $\boldsymbol{x}$ representing the \emph{density}
assigned to one of the pixels in the image. We use $J=235,225$, thus
each $\boldsymbol{x}$ represents a $485\times485$ image. Our test
image (phantom) is represented by the vector $\boldsymbol{\hat{x}}$,
that is a digitization of a picture of a cross-section of a human
head. The picture underlying the digitization is geometrically defined
so that it has a value at every point of the planar cross-section;
see \cite[Sections 4.1--4.4 and 5.2]{GTH-book}.

In the problem $T$ that we use for our illustration, each index $i=1,2,\ldots,I$
is associated with a line across the image and the corresponding $\boldsymbol{d}^{i}$
is a vector in $\mathbb{R}^{J}$, whose $j$th component is the length
of intersection of that line with the $j$th of the $J$ pixels. There
are $I=498,960$ such lines (organized into 720 divergent projections
with 693 lines in each; similar to the \textit{standard geometry}
in \cite{GTH-book} but with more lines in each projection). The $h_{i}$
have been calculated by simulating the behavior of CT scanning of
the head cross-section \cite[Section 4.5]{GTH-book}. In particular,
the projection data were simulated using the underlying picture (rather
than its digitization) and incorporate the stochastic nature (noisiness)
of data collection in CT (based on 1,000,000 photons for estimating
every line integral). All the above was generated using the SNARK14
programming system for the reconstruction of 2D images from 1D projections
\cite{SNARK14}, giving rise to a system of linear equations \eqref{eq:equalities}.
For the resulting $T$, we calculated that the proximity of the phantom
to the generated constraints is $\mathcal{P}r_{T}(\boldsymbol{\boldsymbol{\hat{x}}})=6.4192$,
which is not zero due to the phantom being a digitization of the underlying
picture and the noise incorporated into the calculation the line integrals
$h_{i}$.

For our demonstration we make the simplest choice for $\Delta$, namely,
$\Delta=\Omega=\mathbb{R}^{J}$. Our choice of the target function
$\phi$ is as follows. We index the pixels (i.e., the components of
a vector $\boldsymbol{x}$) by $j$ and let $\Theta$ denote the set
of all indices of pixels that are not in the rightmost column or the
bottom row of the 2D pixel array that displays that vector as an image.
For any pixel with index $j\in\Theta$, let $r\left(j\right)$ and
$b\left(j\right)$ be the index of the pixel to its right and below
it in the 2D pixel array, respectively. Denoting by $\mathrm{med}$
the function that selects the median value of its three arguments,
we define 
\begin{equation}
\phi\left(\boldsymbol{x}\right):=\sum_{j\in\Theta}\sqrt{\left|x_{j}-\mathrm{med}\left\{ x_{j},x_{r\left(j\right)},x_{b\left(j\right)}\right\} \right|}.\label{eq:median_target}
\end{equation}
This function can be considered as an alternative to total variation
for measuring the roughness in an image represented by the vector
$\boldsymbol{x}$. It has been selected for the experiments demonstrating
the behavior of derivative-free superiorization because finding partial
derivatives for it is problematic. On the other hand, when only one
pixel value (that is, only one component of the vector) is changed
in vector $\boldsymbol{x}$ to get another vector $\boldsymbol{y}$,
then it is possible to obtain $\phi(\boldsymbol{y})$ from $\phi(\boldsymbol{x})$
by computing only three of the terms in the summation on the right-hand
side of \eqref{eq:median_target}. These observations indicate that
the use of the derivative-free approach of Steps 10\textendash 18
in \textbf{Algorithm \ref{alg:Component-wise}} is a viable option
whereas Step 8 of \textbf{Algorithm \ref{alg:NonAsc_Superiorization}}
is hard to perform unless the trivial nonascending vector $\boldsymbol{v}^{k,n}=\boldsymbol{0}$
is selected, which is ineffective. For our chosen phantom we calculated
$\phi\left(\boldsymbol{\boldsymbol{\hat{x}}}\right)=2,048.57$.

\subsection{The algorithmic operator $\mathbf{P}_{T}$}

Our chosen operator, mapping $\boldsymbol{x}$ into $\mathbf{P}_{T}\boldsymbol{x}$,
is specified by \textbf{Algorithm \ref{alg:OperatorPT}.} It depends
on a real parameter $\lambda\in(0,2)$ in its Step 4.

\setcounter{algorithm}{2}

\begin{algorithm}

\caption{The algorithmic operator $\mathbf{P}_{T}$}

\label{alg:OperatorPT}

\begin{algorithmic}[1]

\STATE{\textbf{set} $i=0$}

\STATE{\textbf{set} $y^{i}=\boldsymbol{x}$}

\STATE{\textbf{while} $i<I$}

\STATE{$\quad$\textbf{set} $\boldsymbol{y}^{i+1}=\boldsymbol{y}^{i}-\lambda{\displaystyle \frac{\left\langle \boldsymbol{d}^{i},\boldsymbol{y}^{i}\right\rangle -h_{i}}{\Vert\boldsymbol{d}^{i}\Vert^{2}}\boldsymbol{d}^{i}}$}

\STATE{$\quad$\textbf{set} $i=i+1$}

\STATE{\textbf{set} $\mathbf{P}_{T}\boldsymbol{x}=\boldsymbol{y}^{I}$}

\end{algorithmic}

\end{algorithm}

\medskip{}

\begin{algorithm}

\caption{ART (as used in this paper)}

\label{alg:ART}

\begin{algorithmic}[1]

\STATE{\textbf{set} $k=0$}

\STATE{\textbf{set} $\boldsymbol{x}^{k}=\boldsymbol{\bar{x}}$}

\STATE{\textbf{repeat}}

\STATE{$\quad$\textbf{set} $\boldsymbol{x}^{k+1}=\mathbf{P}_{T}\boldsymbol{x}^{k}$}

\STATE{$\quad$\textbf{set} $k=k+1$}

\end{algorithmic}

\end{algorithm}

\textbf{Algorithm \ref{alg:ART}} is a special case of the general
class of Algebraic Reconstruction Techniques as discussed in \cite[Chapter 11]{GTH-book}
and is, for $\lambda=1$, equivalent to the original method of Kaczmarz
in the seminal paper \cite{kacmarcz}. For further references on Kaczmarz's
method and the Algebraic Reconstruction Techniques see, e.g., \cite[page 220]{cegielski-book},
\cite[Section 2]{annotated15} and \cite{HERM19}. Note that \textbf{Algorithm
\ref{alg:ART}} (ART) can be obtained from either \textbf{Algorithm
\ref{alg:NonAsc_Superiorization}} or \textbf{Algorithm \ref{alg:Component-wise}}
by removing the perturbation steps in their \textbf{while} loops.

\subsection{\label{subsec:EP-approach}A comment about the exterior penalty function
approach to derivative-free constrained minimization}

One possibility for doing a derivative-free constrained minimization
algorithm is to follow the option of using the exterior penalty (EP)
function approach mentioned in \cite[Chapter 13, Section 13.1, page 242]{Conn-book-2009}
as applied to the constrained problem 
\begin{equation}
\min\text{ }\text{ }\left\{ \phi\left(\boldsymbol{x}\right)\mid\left\langle \boldsymbol{d}^{i},\boldsymbol{x}\right\rangle =h_{i},\thinspace\thinspace\thinspace i=1,2,\ldots,I\right\} ,\label{eq:constrained}
\end{equation}
where $\phi$ is as in \eqref{eq:median_target} and the constraints
are as in \eqref{eq:equalities}. With a user-selected \emph{penalization
parameter} $\eta$, the exterior penalty function approach replaces
the constrained minimization problem \eqref{eq:constrained} by the
penalized unconstrained minimization:
\begin{equation}
\min\text{ }\{\psi(\boldsymbol{x})\mid\boldsymbol{x}\in\mathbb{R}^{J}\},\label{eq:reg-prob}
\end{equation}
\begin{equation}
\psi(\boldsymbol{x}):=\phi\left(\boldsymbol{x}\right)+\eta\mathcal{P}r_{T}(\boldsymbol{x}),\label{eq:unconstrained}
\end{equation}
with $\phi$ as in \eqref{eq:median_target} and $\mathcal{P}r_{T}(\boldsymbol{x})$
as defined in \eqref{eq:residual}. By applying the coordinate-search
method of \cite[Algorithm 7.1]{Conn-book-2009} to the penalized unconstrained
minimization problem \eqref{eq:reg-prob}\textendash\eqref{eq:unconstrained},
we get the next algorithm.

\begin{algorithm}

\caption{Derivative-free constrained minimization by the exterior
penalty (EP) approach}

\label{alg:Derivative-free}

\begin{algorithmic}[1]

\STATE{\textbf{set} $k=0$}

\STATE{\textbf{set} $\boldsymbol{x}^{k}=\boldsymbol{\bar{x}}$}

\STATE{\textbf{set} $\ell=-1$}

\STATE{\textbf{set} $m=-1$}

\STATE{\textbf{repeat}}

\STATE{$\quad$\textbf{set} $\boldsymbol{x}^{k+1}=\boldsymbol{x}^{k}$}

\STATE{$\quad$\textbf{set} $\ell=\ell+1$}

\STATE{$\quad$\textbf{set} $L=-1$}

\STATE{$\quad$\textbf{while} $L<2J$}

\STATE{$\quad\quad$\textbf{set} $L=L+1$}

\STATE{$\quad\quad$\textbf{set} $m=m+1$}

\STATE{$\quad\quad$\textbf{set} $\boldsymbol{z}=\boldsymbol{x}^{k}+\gamma_{\ell}\mathbf{c}^{m}$}

\STATE{$\quad\quad$\textbf{if} $\boldsymbol{z}\in\Omega$ \textbf{and}
$\psi\left(\boldsymbol{z}\right)<\psi\left(\boldsymbol{x}^{k}\right)$
\textbf{then}}

\STATE{$\quad\quad\quad$\textbf{set} $\boldsymbol{x}^{k+1}=\boldsymbol{z}$}

\STATE{$\quad\quad\quad$\textbf{set} $L=2J$}

\STATE{$\quad$\textbf{set} $k=k+1$}

\end{algorithmic}

\end{algorithm}

From the point of view of keeping the computational cost low, \textbf{Algorithm
\ref{alg:Derivative-free}} can be much more of a challenge than \textbf{Algorithm
\ref{alg:Component-wise}}. The reason for this has been indicated
when we have stated, near the end of Subsection \ref{subsec:The-selected-problem},
that if only one component is changed in vector $\boldsymbol{x}$
to get another vector $\boldsymbol{y}$, then it is possible to obtain
$\phi(\boldsymbol{y})$ from $\phi(\boldsymbol{x})$ by computing
only three of the terms in the summation on the right-hand side of
\eqref{eq:median_target}. When we use $\psi$ in \eqref{eq:unconstrained}
instead of $\phi$, there seems to be a need for many more computational
steps. This is because the number of terms that change on the right-hand
side of \eqref{eq:residual} due to a change in one component of $\boldsymbol{x}$
is of the order of 1,000 for the dataset described in Subsection \ref{subsec:The-selected-problem}
(in the language of image reconstruction from projections, there is
at least one line $i$ in each of of the 720 projections for which
there is a change in value of $\left\langle \boldsymbol{d}^{i},\boldsymbol{x}\right\rangle $
due to changing one component of $\boldsymbol{x}$). Thus our advocated
approach of component-wise superiorization in \textbf{Algorithm \ref{alg:Component-wise}}
is likely to be orders of magnitude faster for our application area
than the more traditional approach of derivative-free constrained
minimization by the exterior penalty (EP) approach in \textbf{Algorithm
\ref{alg:Derivative-free}.}

\subsection{\label{subsec:The-progressive-barrier}The progressive barrier (PB)
approach}

The progressive barrier (PB) approach, described in Chapter 12 of
\cite{Audet-book-2017}, was originally published in \cite{Audet-Dennis-2009}
wherein the history of the approach as a development of the earlier
filter methods of Fletcher and Leyffer \cite{filter02} is succinctly
described. It is an alternative to the exterior penalty (EP) approach,
mentioned above, therefore, we briefly describe it here and point
out how it differs from our DFS. The PB approach bears some similarities
to our DFS but differs from it in a way that explains why the DFS
will be advantageous for large-scale problems. 

No penalty is used in the PB approach. Instead of combining the constraints
with the target function it uses a particular ``constraint violation
function'' $h(x)$ \cite[Definition 12.1]{Audet-book-2017} alongside
with the target function $\phi(x)$ so that the iterates appear in
an $h$ versus $\phi$ plot called ``a filter'', based on the pairs
of their $h$ and $\phi$ values. The constraint violation function
of PB is similar in nature to our ``proximity function'' and the
$h$ versus $\phi$ filter plot of PB is similar to our proximity-target
curve, both mentioned above. The difference between the PB approach
of \cite{Audet-Dennis-2009} and our DFS lies in how these objects
are used. The PB optimization algorithm defines at each iteration
what is a ``success'' or a ``failure'' of an iterate based on
the current filter plot, and decides accordingly what will be the
next iterate.  We do not bring the full details here but, in a nutshell,
the PB optimization algorithm performs at each iteration sophisticated
searches of both the target function values and the constraint violation
function values. 

In contrast with the PB approach, the DFS, investigated here, uses
the world-view of superiorization. It uses a feasibility-seeking algorithm
whose properties are already known and perturbs its iterates without
losing its feasibility-seeking ability and properties. The component-wise
derivative free search for a locally nonascending direction of the
target function is done in a manner that makes it a perturbation to
which the feasibility-seeking algorithm is resilient, i.e., that allows
the underlying feasibility-seeking algorithm retain its feasibility-seeking
nature. Thus, the DFS method searches only the target function in
a derivative-free fashion and does automatic feasibility-seeking steps
that actively reduce the proximity function. In this work the constraints
are linear and the feasibility-seeking algorithm proceeds by performing
orthogonal projections onto the hyperplanes that constitute the constraints
sets). This means, on the face of it, an advantage for the DFS in
handling large-scale problems because no expensive additional time
and computing resources are needed for the feasibility-seeking phase
of the DFS algorithm proposed here. The validity of this point requires
further research.

\subsection{\label{subsec:Algorithmic-details-and}Algorithmic details and numerical
demonstration}

Our experiments were carried out using the public-domain software
package SNARK14 \cite{SNARK14}. In all experiments the initial vector
$\boldsymbol{\bar{x}}$ was the $235,225$-dimensional zero vector
(all components 0).

The relaxation parameter in \textbf{Algorithm \ref{alg:OperatorPT}}
was $\lambda=0.05$. Another issue that needs specification is the
ordering of the constraints in \eqref{eq:equalities}, because the
output of \textbf{Algorithm \ref{alg:OperatorPT}} depends not only
on the set of constraints, but also on their order. We used in our
experiments the so-called \emph{efficient ordering}, since it has
been demonstrated to lead to better results faster when incorporated
into ART \cite[page 209]{GTH-book}.

\begin{figure}[t]
\begin{centering}
\includegraphics[width=1\columnwidth]{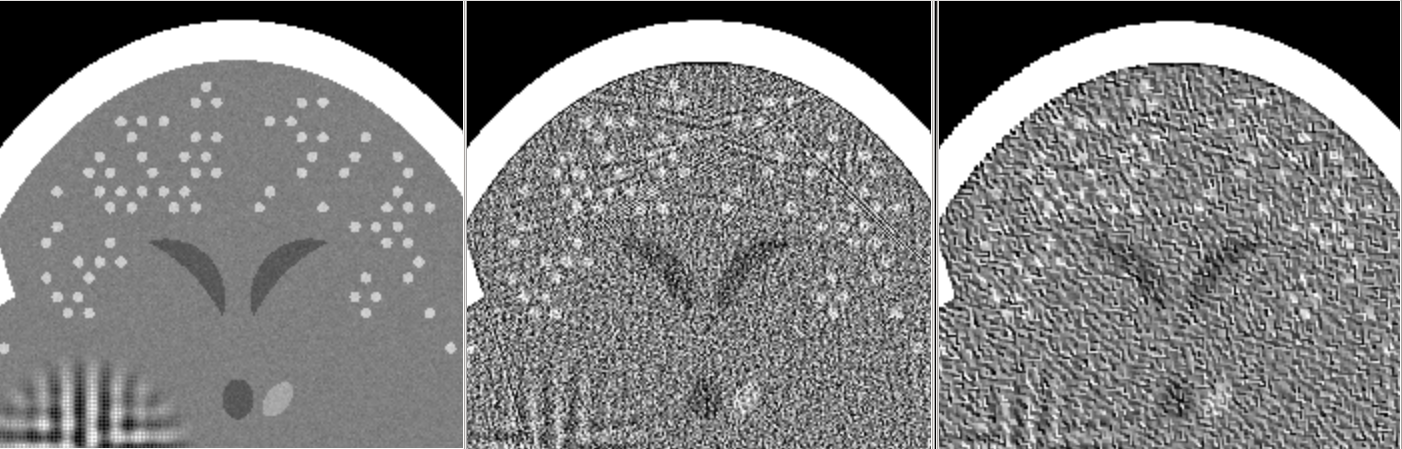}
\par\end{centering}
\caption{\label{fig:Details_Exp}Detail of (left) a 485\texttimes 485 digitization
of a phantom based on the distribution of x-ray attenuation in units
of cm$^{-1}$ within a transaxial slice of the human head (center)
reconstruction produced by \textbf{Algorithm \ref{alg:Component-wise}}
without its component-wise perturbations steps, and (right) produced
by \textbf{Algorithm \ref{alg:Component-wise}} with its component-wise
perturbations steps. In all images we display a value that is 0.20
or less as black and a value that is 0.22 or greater as white.}
\end{figure}

In \textbf{Algorithm \ref{alg:Component-wise}}, the number $N$ of
perturbation steps (for each feasibility-seeking step) was 100,000
and we used $\gamma_{\ell}=ba^{\ell}$, with $b=0.02$ and $a=0.999,999$.
The infinite sequence $\left(\mathbf{c}^{m}\right)_{m=0}^{\infty}$
was obtained by repetitions of the length-$2J$ subsequence $\left(\mathbf{e}^{1},\mathbf{e}^{2},\ldots\mathbf{e}^{J},-\mathbf{e}^{1},\mathbf{-e}^{2},\ldots,-\mathbf{e}^{J}\right)$.

\begin{figure}
\begin{centering}
\includegraphics[scale=0.15]{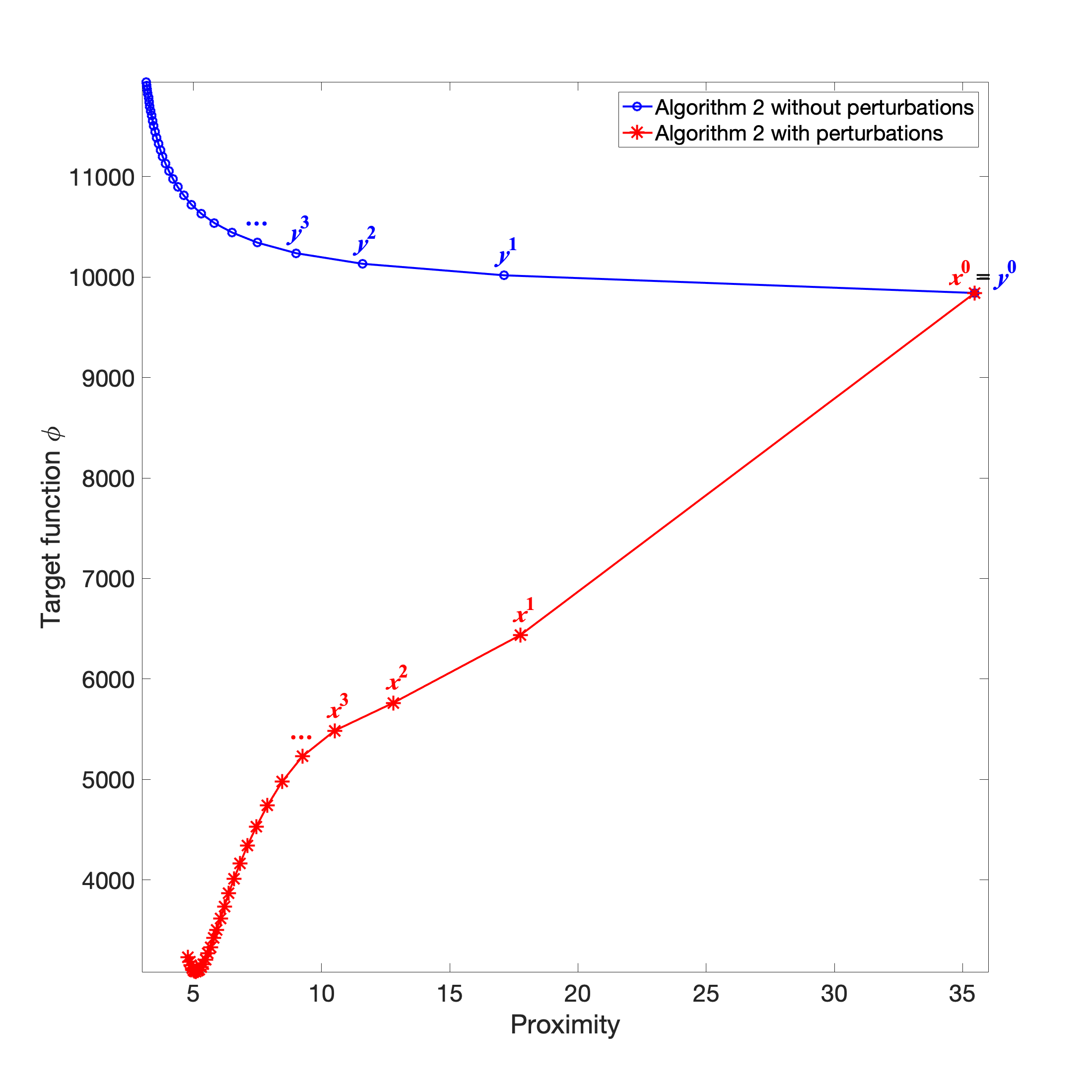}
\par\end{centering}
\caption{\label{fig:Target-function-values}Proximity-target curves $P$ and
$Q$ of the first 30 iterates of \textbf{Algorithm 2} with perturbations
(\EightStarTaper ) and without perturbations ($\circ)$.}
\end{figure}

We applied \textbf{Algorithm \ref{alg:Component-wise}} twice, thirty
iterations in each case, with and without its component-wise perturbations
steps, respectively, under otherwise completely identical conditions;
in Figure \ref{fig:Details_Exp} we show a detail of the results produced
by these two executions together with the corresponding detail of
the digitized phantom. The resulting finite sequences of iterates
are both of monotone proximity, the associated proximity-target curves
are shown in Figure \ref{fig:Target-function-values}. The $\circ$s
and \EightStarTaper s on the plots represent actually calculated values
at iterations of each algorithm, that are connected by line segments.
For any proximity value on the horizontal axis we can read the target-function
value associated with it from the curve. The plots indicate visually
the behavior of the algorithms, initialized at the same point denoted
by $x^{0}=y^{0}$ that appears in the right-most side of the figure.
The V-shaped form of the proximity-target curve for \textbf{Algorithm
\ref{alg:Component-wise}} with perturbations is typical for the behavior
of superiorized feasibility-seeking algorithms, showing the initially
strong effect of the perturbations that diminishes as the iterations
proceed.

For this experiment we used SNARK14 \cite{SNARK14} installed in a
computer with a 2.7 GHz Intel\textsuperscript{\textregistered} Core
64-bit i7 processor and a total main memory of 16 Gbytes running the
CentOS Linux operating system on a virtual machine with 8 Gbytes assigned
for base memory. Under these conditions, the implementation of \textbf{Algorithm
\ref{alg:Component-wise}} without component-wise perturbations steps
took 77.251 seconds to execute all the iterations whereas the implementation
with such perturbations took 105.801 seconds for all the iterations.
The latter time was measured excluding the time the implementation
took to compute the useless Steps 12 to 18 for the first iteration
(i.e., $k=0$). These steps are useless in our case, since with the
zero vector as initial value, the condition in Step 16 is never satisfied
while $k=0$, resulting in $\boldsymbol{x}^{0,N}=\boldsymbol{x}^{0}$.

For a more precise interpretation, consider Definition \ref{def:comparison}.
In the experiment evaluating the two versions of \textbf{Algorithm
\ref{alg:Component-wise}}, $K_{lo}=L_{lo}=1$ and $K_{hi}=L_{hi}=30$.
The $R=\left(\boldsymbol{x}^{k}\right)_{k=K_{lo}}^{K_{hi}}$ and $S=\left(\boldsymbol{y}^{k}\right)_{k=L_{lo}}^{L_{hi}}$
produced by \textbf{Algorithm \ref{alg:Component-wise}}, with and
without perturbations, respectively, are both of monotone proximity.
We find that $\mathcal{P}r_{T}\left(\boldsymbol{x}^{K_{lo}}\right)=Pr_{T}\left(\boldsymbol{y}^{L_{lo}}\right)=35.4703$
(and, hence, $u=35.4703$) and that $\mathcal{P}r_{T}\left(\boldsymbol{x}^{K_{hi}}\right)=3.4065$
and $\mathcal{P}r_{T}\left(\boldsymbol{y}^{L_{hi}}\right)=4.7828$
(and, hence, $t=4.7828$). By showing the target curves $P$ and $Q$
associated with $R$ and $S$, respectively, Figure \ref{fig:Target-function-values}
clearly illustrates that $R$ is better targeted than $S$.

\section{\label{sec:Discussion-and-conclusions}Discussion and conclusions}

In this paper we investigated the general applicability of derivative-free
superiorization (DFS) as an alternative to previously proposed superiorization
approaches. In our computational demonstration, we generated the constraints
and chose the target function from the application area of image reconstruction
from projections in computerized tomography (CT). However, we use
the demonstration for indicating only the numerical behavior of the
algorithms. We do not investigate or comment on the potential usefulness
of the resulting reconstructions in CT, since that usefulness depends
not so much on the numerical behavior of the algorithms as on the
appropriateness of the modeling used to turn a physical problem into
a mathematical one (for example, by the specific choice of target
function). The numerical results of our demonstration attest, as seen
from the proximity-target curves, to the mathematical efficacy of
our derivative-free superiorization algorithm, but say nothing about
its efficacy for providing an answer to a practical image reconstruction
problem. (Nevertheless, we have observed while doing our experiment
that, even from the image reconstruction quality point of view, DFS
seems to be advantageous. For example, if we consider the distances
between the phantom and the reconstructions -defined as the 2-norm
between the representing vectors-, the smallest distance that we get
as we iterate without perturbations is 0.0922, while with the DFS
perturbations it is 0.0863.)

Much of the literature on derivative-free minimization is concerned
with unconstrained minimization or at most with bound-constraints
on the variables, and many, if not all, proposed methods can handle
only small-size problems efficiently. In contrast, the DFS method
proposed here can handle any type of constraints for which a separate
efficient derivative-free constraints-compatibility-seeking algorithm
is available. Since the constraints-compatibility-seeking algorithm
forms part of the proposed DFS method, the method can use exterior
initialization (i.e., initializing the iterations at any point in
space). Furthermore, and very importantly, very large-size problems
can be accommodated.
\begin{acknowledgements}
We thank Nikolaos Sahinidis and Katya Scheinberg for several informative
mail exchanges that helped us see better the general picture. We are
grateful to Sébastien Le Digabel for calling our attention to the
work of Charles Audet and coworkers, particularly the Audet and Dennis
paper \cite{Audet-Dennis-2009} and the book of Audet and Hare \cite{Audet-book-2017}.
We greatly appreciate the constructive referee report that helped
us improve the paper.
\end{acknowledgements}

\section*{Conflict of interest}

The authors declare that they have no conflict of interest.

\bibliographystyle{spmpsci}

\begin{thebibliography}{10}
\bibitem{Audet-Dennis-2009}Audet, C. and Dennis J.E. JR., 2009. A
progressive barrier for derivative-free nonlinear programming, \textit{SIAM
Journal on Optimization}, 20, 445\textendash 472.

\bibitem{Audet-book-2017}Audet, C. and Hare, W., 2017. \textit{Derivative-Free
and Blackbox Optimization}, Springer International Publishing, Cham,
Switzerland.

\bibitem{Bargetz2018}Bargetz, C., Reich, S., and Zalas, R., 2018.
Convergence properties of dynamic string-averaging projection methods
in the presence of perturbations, \textit{Numerical Algorithms}, 77,
185\textendash 209.

\bibitem{But06} Butnariu, D., Reich, S., and Zaslavski, A. J., 2006.
Convergence to fixed points of inexact orbits of Bregman-monotone
and of nonexpansive operators in Banach spaces, \textit{in: Proceedings
of Fixed Point Theory and its Applications}\emph{, Mexico}, Yokohama,
11\textendash 32.

\bibitem{cegielski-book}Cegielski, A., 2012. \textit{Iterative Methods
for Fixed Point Problems in Hilbert Spaces}, Springer-Verlag.

\bibitem{Cegielski-2017}Cegielski, A. and Al-Musallam, F., 2017.
Superiorization with level control, \textit{Inverse Problems}, 33,
044009.

\bibitem{Cen15}Censor, Y., 2015. Weak and strong superiorization:
Between feasibility-seeking and minimization, \textit{Analele Stiintifice
ale Universitatii Ovidius Constanta-Seria Matematica}, 23, 41\textendash 54.

\bibitem{sup-bib}Censor, Y., 2019. \textit{Superiorization and Perturbation
Resilience of Algorithms: A Bibliography compiled and continuously
updated}, \href{http://math.haifa.ac.il/yair/bib-superiorization-censor.html}{http://math.haifa.ac.il/yair/bib-superiorization-censor.html},
last updated: September 19, 2020.

\bibitem{annotated15}Censor, Y. and Cegielski, A., 2015. Projection
methods: An annotated bibliography of books and reviews, \textit{Optimization},
64, 2343\textendash 2358.

\bibitem{CHS18}Censor, Y., Heaton, H., and Schulte, R.W., 2019. Derivative-free
superiorization with component-wise perturbations, \textit{Numerical
Algorithms}, 80, 1219\textendash 1240.

\bibitem{Sup-Special-Issue-2017}Censor, Y., Herman, G.T., and Jiang,
M., (Editors) 2017. Special issue on Superiorization: Theory and Applications,
\textit{Inverse Problems}, 33, 040301\textendash 044014.

\bibitem{censor-levy-2019}Censor, Y., Levy, E., 2019. An analysis
of the superiorization method via the principle of concentration of
measure, \textit{Applied Mathematics and Optimization}. https://doi.org/10.1007/s00245-019-09628-4.

\bibitem{cz3-2015}Censor, Y. and Zaslavski, A., 2015. Strict Fejér
monotonicity by superiorization of feasibility-seeking projection
methods, \textit{Journal of Optimization Theory and Applications},
165, 172\textendash 187.

\bibitem{Conn-book-2009}Conn, A.R., Scheinberg, K., and Vicente,
L.N., 2009. \textit{Introduction to Derivative-Free Optimization},
Society for Industrial and Applied Mathematics (SIAM).

\bibitem{COOL92a}Cooley, T.A. and Barrett, H.H., 1992. Evaluation
of statistical methods for image reconstruction through ROC analysis,\textit{
IEEE Transactions on Medical Imaging}, 11, 276-282.

\bibitem{SNARK14}Davidi, R., Garduño, E., Herman, G.T., Langthaler,
O., Rowland, S.W., Sardana, S., and Ye, Z., 2019. SNARK14: A programming
system for the reconstruction of 2D images from 1D projections, available
from \href{http://turing.iimas.unam.mx/SNARK14M/SNARK14.pdf}{http://turing.iimas.unam.mx/SNARK14M/SNARK14.pdf}.

\bibitem{diniz2011}Diniz-Ehrhardt, M., Martínez, J., and Pedroso,
L., 2011. Derivative-free methods for nonlinear programming with general
lower-level constraints, \textit{Computational and Applied Mathematics},
30, 19\textendash 52.

\bibitem{dfo-4-oil}Echeverría Ciaurri, D., Isebor, O., and Durlofsky,
L., 2012. Application of derivative-free methodologies to generally
constrained oil production optimization problems, \textit{Procedia
Computer Science}, 1, 1301\textendash 1310.

\bibitem{ehn00}Engl, H.W., Hanke, M., and Neubauer, A., 2000. \textit{Regularization
of Inverse Problems}, Kluwer Academic Publishers.

\bibitem{filter02}Fletcher, R. and Leyffer, S., 2002. Nonlinear programming
without a penalty function, \textit{Mathematical Programming, Series
A}, 91, 239\textendash 269.

\bibitem{GH14}Garduño, E. and Herman, G.T., 2014. Superiorization
of the ML-EM algorithm, \textit{IEEE Transactions on Nuclear Science},
61, 162\textendash 172.

\bibitem{GTH}Garduño, E. and Herman, G.T., 2017. Computerized tomography
with total variation and with shearlets, \textit{Inverse Problems},
33, 044011.

\bibitem{gay-niemierko-ntcp-2007}Gay, H.A., Niemierko, A., 2007.
A free program for calculating EUD-based NTCP and TCP in external
beam radiotherapy, \textit{Physica Medica}, 23, 115\textendash 25.

\bibitem{He2017}He, H. and Xu, H.K., 2017. Perturbation resilience
and superiorization methodology of averaged mappings, \textit{Inverse
Problems}, 33, 044007.

\bibitem{GTH-book}Herman, G.T., 2009. \textit{Fundamentals of Computerized
Tomography: Image Reconstruction from Projections}, Springer-Verlag,
2nd ed.

\bibitem{HERM19}Herman, G.T., 2019. Iterative reconstruction techniques
and their superiorization for the inversion of the Radon transform,
\textit{in}: R. Ramlau and O. Scherzer, eds., \textit{The Radon Transform:
The First 100 Years and Beyond}, De Gruyter, 217\textendash 238.

\bibitem{Herman-JANO}Herman, G.T., 2020. Problem structures in the
theory and practice of superiorization, \textit{Journal of Applied
and Numerical Optimization}, 2, 71\textendash 76.

\bibitem{HGDC12}Herman, G.T., Garduño, E., Davidi, R., and Censor,
Y., 2012. Superiorization: An optimization heuristic for medical physics,
\textit{Medical Physics}, 39, 5532\textendash 5546.

\bibitem{hoseini-2019}Hoseini, M., Saeidi, S. and Kim, D.S., 2019.
On perturbed hybrid steepest descent method with minimization or superiorization
for subdifferentiable functions. \textit{Numerical Algorithms}. https://doi.org/10.1007/s11075-019-00818-3.

\bibitem{kacmarcz}Kaczmarz, S., 1937. Angenäherte Auflösung von Systemen
linearer Gleichungen, \textit{Bulletin de l'Académie Polonaise des
Sciences et Lettres}, A35, 355\textendash 357.

\bibitem{Kolda2003}Kolda, T., R. Lewis and V. Torczon, 2003. Optimization
by direct search: New perspectives on some classical and modern methods,
\textit{SIAM Review,} 45, 385\textendash 482.

\bibitem{LCLLLL}Li, L., Chen, Y., Liu, Q., Lazic, J., Luo, W., and
Li, Y., 2017. Benchmarking and evaluating MATLAB derivative-free optimisers
for single-objective applications, \textit{in}: D.S. Huang, K.H. Jo,
and J. Figueroa-García, eds., \textit{Intelligent Computing Theories
and Application. ICIC 2017}, Springer, 75\textendash 88.

\bibitem{LZZS16}Luo, S., Zhang, Y., Zhou, T., and Song, J., 2016.
Superiorized iteration based on proximal point method and its application
to XCT image reconstruction, \textit{ArXiv e-prints}, \href{https://arxiv.org/abs/1608.03931}{https://arxiv.org/abs/1608.03931}.

\bibitem{LZZSW19}Luo, S., Zhang, Y., Zhou, T., Song, J., and Wang,
Y., 2019. XCT image reconstruction by a modified superiorized iteration
and theoretical analysis, \textit{Optimization Methods and Software},
DOI: 10.1080/10556788.2018.1560442.

\bibitem{METZ89a}Metz, C.E., Some practical issues of experimental
design and data analysis in radiological ROC studies, \emph{Investigative
Radiology}, 24, 234-243.

\bibitem{NDH12}Nikazad, T., Davidi, R., and Herman, G.T., 2012. Accelerated
perturbation resilient block-iterative projection methods with application
to image reconstruction, \textit{Inverse Problems}, 28, 035005.

\bibitem{nystorm-2020}Nystrom, H., Jensen, M.F., Nystrom, P.W., 2020.
Treatment planning for proton therapy: what is needed in the next
10 years? \textit{The British Journal of Radiology}, 93, 1107, 20190304.
DOI:10.1259/bjr.20190304.

\bibitem{reich-zalas-2016}Reich, S. and Zalas, R., 2016. A modular
string averaging procedure for solving the common fixed point problem
for quasi-nonexpansive mappings in Hilbert space, \textit{Numerical
Algorithms}, 72, 297\textendash 323.

\bibitem{Reich2017}Reich, S. and Zaslavski, A.J., 2017. Convergence
to approximate solutions and perturbation resilience of iterative
algorithms, \textit{Inverse Problems}, 33, 044005.

\bibitem{Rios}Rios, L.M. and Sahinidis, N.V., 2013. Derivative-free
optimization: A review of algorithms and comparison of software implementations,
\textit{Journal of Global Optimization}, 56, 1247\textendash 1293.

\bibitem{SWET79a}Swets, J.A.,1979. ROC analysis applied to the evaluation
of medical imaging techniques, \emph{Investigative Radiology}, 14,
109\textendash 112.

\bibitem{ZAS18}Zaslavski, A.J., 2018. \emph{Algorithms for Solving
Common Fixed Point Problems}, Springer International Publishing.
\end{thebibliography}

\addcontentsline{toc}{section}{\refname}
\end{document}